\documentclass[11pt,leqno]{article}
\usepackage[latin1]{inputenc}
\usepackage{amsxtra}
\usepackage{amsmath}
\usepackage{amssymb}
\usepackage{amsfonts}
\usepackage[all]{xy}
\usepackage{mathrsfs}
\usepackage{amsthm}
\usepackage{enumerate}

\newtheorem{lem}[equation]{Lemma}
\newtheorem{prop}[equation]{Proposition}

\newtheoremstyle{example}{\topsep}{\topsep}%
     {}
     {}
     {\bfseries}
     {.}
     {2pt}
     {\thmname{#1}\thmnumber{ #2}\thmnote{ #3}}

   \theoremstyle{example}

   \newtheorem{ex}[equation]{Example}

  \numberwithin{equation}{section}

\def\CC{\mathbb{C}}

\def\PP{\mathbb{P}}
\def\RR{\mathbb{R}}
\def\ZZ{\mathbb{Z}}

 \def\Conc{\operatorname{Conc}}

\def\3x3{\operatorname{3}\times\operatorname{3}}
\def\1{{\bf 1}}
\def\lra{\longrightarrow}

\def\(({(\hskip -1mm (}
\def\)){)\hskip -1mm )}
\def\be{\begin{equation}}
\def\ee{\end{equation}}
\def\ed{\end{document}}

\def\la{{\langle}}
\def\ra{{\rangle}}

\title{ Thermodynamics and the moment map}

\author{ M. Kapranov}

\begin{document}

\maketitle

\begin{abstract}
We give a thermodynamic interpretation of the moment map for
toric varieties. The convexity properties of this map correspond
to  thermodynamical principles (concavity of the entropy functional)
applied to a system with several Hamiltonians. 
\end{abstract}


\addtocounter{section}{-1}

\section {Introduction.}

This elementary note is an exercise in  generalizing
the standard (Maxwell-Bolztmann-Gibbs) approach to thermodynamics,
to the case when the energy function is vector valued. This case is
similar to that of several commuting Hamiltonians, familiar
in the theory of integrable systems. 

It turns out that the mean energy function of such ``higher-dimensional
thermodymanics" is basically the same as the moment map in the theory
of toric varieties. Low-temperature limit of the standard theory
generalizes to the ``tropical limits" near vertices of the convex polytope
given by the image of the moment map.  Convexity properties of
the moment map correspond to
 fundamental thermodynamic principles
such as concavity of the entropy functional.  

Relation of thermodynamics with tropical geometry have recently begun to
attract some interest  from various directions \cite{itenberg-mikhalkin, marcolli}. In particular, 
the observation that tropical geometry corresponds, thermodynamically, to the 
low-temperature limit (in contrast with the name ``tropical"
which suggests the opposite)
has been made  by I. Itenberg and G. Mikhalkin in \cite{itenberg-mikhalkin}.
Consideration of vector inverse temperature as in this note, 
 makes this observation even more clear.

I am gratefiul to G. Mikhalkin for stimulating discussions.

\section{Reminder on usual thermodynamics (one Hamiltonian).}
We start by reviewing the standard material on statistical thermodynamics,
see \cite{landau-lifschitz, zinsmeister}, with some extra emphasis on logical structure. 

\paragraph{A. The Gibbs distribution.}
Consider a thermodynamical system such as a gas, with the (large but finite) set of states $A$.
According to the fundamental principle of Boltzmann, the statistical behavior of the system
is completely determined by the knowledge of the energies of various states, i.e., by
the choice of a function $E: A\to\RR$. 
%
%
To find this behavior, one forms the {\em Boltzmann partition function}
\be\label{eq:partition-one}
Z(\beta) \,\,=\sum_{\omega\in A} e^{-\beta\cdot E(\omega)},\quad \beta={1\over kT},
\ee
where $T$ is the absolute temperature and $k$ is the Boltzmann constant. 
This is a finite sum of exponents, so it is  well-defined    and positive for all real values
of $\beta$. The {\em Gibbs probability distribution} is given by
\be\label {eq:gibbs-one}
p_\omega(\beta) \,\,=\,\,{e^{-\beta\cdot E(\omega)}\over Z(\beta)}, \quad \quad
\sum_\omega p_\omega(\beta)=1.
\ee
The interpretation of  $p_\omega(\beta)$ is  as follows. Let us heat the system
to  temperature $T=1/k\beta$ and wait till it arrives   at a ``thermodynamical equilibrium".
Then $p_\omega(\beta)$ is the probability that the system is at the state $\omega$. 

By an {\em observable} we mean simply a function $O: A\to\RR$. The Gibbs distribution
gives rise to the {\em mean value} of $O$, which is a function
\be\label{eq:mean-value-O}
\la O\ra\,=\,\la O\ra (\beta)\,\,=\,\,\sum_{\omega\in A} p_\omega(\beta)\cdot O(\omega).
\ee
In particular, we have the mean value of the energy
\be\label{eq:mean-energy-beta-1}
\la E\ra (\beta) \,\,=\,\, {\sum_{\omega\in A} E(\omega)\cdot e^{-\beta E(\omega)}
\over \sum_{\omega'\in A} e^{-\beta E(\omega')} } \,\,=\,\, -{d\over d\beta} \log \, Z(\beta).
\ee
Let $E_{\min}, E_{\max}$ be the minimal and maximal values of $E$. For
simplicity assume that each of these values is attained at one state: $\omega_{\min}, \omega_{\max}$.
In the low temperature limit $\beta\to +\infty$ the value $\la E\ra (\beta)$ approaches $E_{\min}$,
as $p_{\omega_{\min}}(\beta)$ approaches 1.  The state $\omega_{\min}$  usually has a
clear physical meaning and is called the {\em ground state}. 
Similarly, in the other\footnote{This is not the high temperature
limit but rather the non-physical limit of $T$ approaching $0$ from the negative 
direction. The limit $T\to +\infty$ corresponds to $\beta\to 0$, when all the states become
equally probable.}
 limit $\beta\to -\infty$ we have that
$\la E\ra (\beta)\to E_{\max}$.  In fact, we have

\begin{prop}\label{prop:diffeomorphsm-1}
The function $\la E\ra(\beta)$ defines a monotone decreasing diffeomorphism from $\RR$
to the open interval $(E_{\min}, E_{\max})$. 
\end{prop}

\noindent {\sl Proof:}
 We need only to show that  $\la E\ra (\beta)$ is
monotone decreasing.  But 
$$
\la E\ra'(\beta)\,\,=\,\,{1\over Z(\beta)^2} \sum_{\omega, \omega'\in A} \bigl( -E(\omega)^2 + E(\omega) E(\omega')\bigr)
e^{-\beta(E(\omega)+E(\omega'))} .
$$
For an unordered pair $\{\omega\neq\omega'\}$ the two coefficients at $e^{-\beta(E(\omega)+E(\omega'))}$
sum to $-(E(\omega)-E(\omega'))^2\leq 0$, while for $\omega=\omega'$ the coefficient vanishes. 
So unless all the $E(\omega)$ are equal to each other, the derivative is strictly negative.  \qed

\paragraph{B. Derivation of the Gibbs distribution.} For future convenience,
we sketch here the classical derivation of \eqref{eq:gibbs-one}. As many thermodynamical arguments,
it assumes {\em two levels of microscopicity}. That is, although the set $A$ 
is already supposed to be very large, $|A| \gg 0$, and involve microscopic degrees of freedom,
 we now assume that we have a much larger
number $N\gg |A|$ of  ``truly microscopic"
particles which can be distributed among the states $\omega\in A$,
possibly many at a time. Each such way of distributing particles is called a
{\em microstate} \footnote{
The ``thermostat" in 
standard discussions of equilibrium thermodynamics is a device for
producing these microstates.
}. 
We further assume (this  corresponds to the classical and not quantum
approach to the problem) that the   particles   distributed are distinguishable
from each other.  This means that we can think of microstates as being
sequences $\xi=(\xi_1, ..., \xi_N )$ of elements of $A$, forming
the Cartesian power $A^N$.

\vskip .2cm

Since we want to determine a probability distribution (measure) on $A$, consider first
the space $\Delta^A$ of all such measures. This is a simplex of dimension $|A|-1$. 
Fix an arbitrary $p=(p_\omega)_{\omega\in A}\in\Delta^A$ and equip $A^N$ with the
product measure. For a microstate $\xi\in A^N$ as above let $N_\omega(\xi) \,=\,|\{i: \xi_i=\omega\}|$
be the number of particles in the state $\omega$, so the  {\em observed probability} 
of being in the state $\omega$ (observed at $\xi$) is $q_\omega(\xi)= N_\omega(\xi)/N$. Fix
a  partition $N=\sum_{\omega\in A} n_\omega$, $n_\omega\in \ZZ_+$. Then the
set of $\xi\in A^N$ such that $N_\omega(\xi)=n_\omega$ has the measure (probability)
\be\label{eq:number-xi}
N! \prod_{\omega\in A} {p_\omega^{n_\omega}\over n_\omega!},
\ee
as the measure of any single such $\xi$ is $\prod p_\omega^{n_\omega}$, while the
number of these $\xi$ is the multinomial coefficient. Using the Stirling approximation
\be
\log(n!) \,\,\sim\,\, n(\log(n)-1), \quad n\gg 0,
\ee
we approximate the logarithm of \eqref{eq:number-xi} by
$$N\sum_{\omega\in A} q_\omega(\log(p_\omega)-\log(q_\omega)),\quad q_\omega = n_\omega/N.$$
Notice the following fact.

\begin{lem} Let $p\in\Delta^A$ be given. Then
$$\max_{q\in\Delta^A} \,\, \sum_{\omega\in A} q_\omega(\log(p_\omega)-\log(q_\omega))=0,
$$
and the maximum is achieved for $q=p$. 
\end{lem}

\noindent {\sl Proof:} The function of $q\in\Delta^A$
which we seek to maximize,    is concave, approaches
$-\infty$ at the boundary and has a critical point at $q=p$, which must then be the absolute maximum.
\qed 

\vskip .1cm

Therefore\footnote
{This is, essentially, the law of large numbers of
probability theory. Up to now, consideration of microstates
was formally identical with that  of  $N$ independent
trials of a random event such as a roll of dice.
} 
 the most probable microstates will be those $\xi$ for which each $q_\omega(\xi)
= p_\omega$. Such $\xi$  are called {\em equilibrium microstates}.
In the case when the $p_\omega=n_\omega/N$ are rational, their number
is the multinomial coefficient, which we interpolate for arbitrary $p\in\Delta^A$,
using the Gamma function,  by
\be\label{eq:microstates}
\text{Number of equilibrium microstates} \,\,\sim \,\, 
{\Gamma(N+1)\over \prod_\omega \Gamma (N p_\omega+1)}.
\ee
 Using the Stirling formula, we approximate the logarithm of \eqref{eq:microstates}
 by
  $$
 N(\log(N)-1) -\sum_\omega Np_\omega
(\log(Np_\omega)-1)) \,\,=\,\, -N \sum_\omega p_\omega \log(p_\omega).
$$
Recall that for a probability distribution $p\in\Delta^A$ its {\em entropy} is defined as
\be\label{eq:entropy-def}
S(p) \,\,=\,\,-\sum_{\omega\in A} p_\omega \log(p_\omega).
\ee
The function $S$ is a concave function on the simplex $\Delta^A$, equal to 0
at the boundary,
and achieving the maximim at the  barycenter. Thus, thermodynamically,
\be
S \,\,\approx\,\,  {\log(\text{Number of equilibrium microstates})\over N
}, \quad N\gg |A|.
\ee

Now, the main thermodynamic principle used to deduce the Gibbs distribution
is that {\em the number of equilibrium microstates  should be as large as possible,
while maintaining the desired mean value of energy}. 
That is, take a point $\overline E\in (E_{\min}, E_{\max})$ and look at all
probability distributions
$p\in\Delta^A$ satisfying
 \be\label{eq:energy-constraint-one}
 \la E\ra_p \,\,:= \,\, \sum_\omega p_\omega E(\omega)\,\, =\,\, \overline E. 
 \ee
 The above principle implies \eqref{eq:gibbs-one} in virtue of the following fact.
 
 \begin{prop} Among the distributions $p$ satisfying \eqref{eq:energy-constraint-one},
 the maximal entropy is achieved by the Gibbs distribution $p(\beta)$,
 where $\beta\in\RR$ is the unique
   number such that $\la E\ra(\beta)$ as defined by 
 \eqref{eq:mean-energy-beta-1}, is equal to $\overline E$. 
  \end{prop}
  
  \noindent {\sl Proof:} The constraint  \eqref{eq:energy-constraint-one}
  defines a hyperplane section of the simplex $\Delta^A$, a convex polytope,
  denote it $P$. The restriction $S|_P$ is a strictly concave function,
  equal to 0 at the boundary. So it has a unique critical
  point inside $P$, and this point is the global maximum.  
  By the Lagrange multiplier method, this critical point is characterized as a
  point $p\in\Delta^A$  which satisfies the constraint (lies in $P$)
  and at which the differential of $S$ is proportional to the differential of the
  constraint. On $\Delta^A$ we have $\sum_\omega dp_\omega=0$, therefore
   $dS=-\sum_\omega \log(p_\omega) dp_\omega$, and the condition of
   proportionality reads:
  \be\label{eq:proportionality}
   -\sum_\omega \log(p_\omega) dp_\omega \,\,=\,\,\lambda \sum_\omega E(\omega) dp_\omega.
   \ee
   But this condition is satisfied by $p_\omega=p_\omega(\beta)$ as defined in the statement
   of the proposition, with $\lambda=\beta$. Indeed, $\log p_\omega(\beta) =-\beta E(\omega) -
   \log Z(\beta)$, so in virtue of $\sum_\omega dp_\omega=0$, the LHS of
   \eqref{eq:proportionality} is equal to $\beta\sum E(\omega) dp_\omega$.
   \qed

\paragraph {C. Entropy, energy and temperature.} One can object that
the above derivation of the Gibbs
distribution \eqref{eq:gibbs-one} is somewhat circular.  It does
 explain, from clear  principles, the behavior of $p=(p_\omega)_{\omega\in A}$
as a function of  the mean energy $\overline E$, but not of $\beta$ or of temperature.
In fact,  $\overline E$ and
$\beta$ are supposed to be related by the formula \eqref{eq:mean-energy-beta-1}
which depends on \eqref{eq:gibbs-one}. This is not surprising since we have not
used any meaningful features of the concept of ``temperature". 

A mathematically satisfying way of dealing with this issue  is to consider the temperature
as a secondary quantity, and to {\em define} it in terms of more fundamental quantities
such as energy.  A standard definition like this
 (see,   e.g.,  \cite{landau-lifschitz})
says that the inverse temperature (i.e., $\beta$) is ``the derivative of the entropy 
with respect to
the energy". Mathematically, this definition (or, rather, its consistency)
amounts to the following general fact
about exponential sums.

\begin{prop} Consider $\beta$ as a function of $\overline E\in (E_{\min}, E_{\max})$
by inverting the diffeomorphism of Proposition 
\ref{prop:diffeomorphsm-1}.  Let $S(\overline E)$ be the entropy of the Gibbs
distribution $p(\beta(\overline E))$. Then $dS/d\overline E = \beta(\overline E)$. 
\end{prop}

In particular, $S(\overline E)$ is a concave function on $[E_{\min}, E_{\max}]$,
equal to $0$ at both ends, with the derivative at these points being $\pm\infty$. 

\vskip .2cm

\noindent {\sl Proof:} This is an immediate consequence of 
 the identification of $\beta$
with the Lagrange multiplier in \eqref{eq:proportionality}.
Indeed, for any constrained maximum problem
$\max_{g(x)=c} f(x)$ the value of the Lagrange multiplier $\lambda=\lambda(c)$
 has the interpretation as the derivative, with respect to $c$,  of the 
 constrained maximum value (this derivative is called, in the language
 of applied math, the ``effective price of the resource
 represented by the constraint", 
 see,  e.g., \cite{econ-book}). 
  \qed

 \section{Thermodynamics with several Hamiltonians.}
 
 \paragraph{A. The Gibbs distribution for several Hamiltonians.}
 We now assume that the set of states $A$ is equipped with
 not one, but several ``energy functionals" $E_1, ..., E_n: A\to\RR$,
 which we combine into one vector valued function $E: A\to\RR^n$. 
 To these energy functionals there correspond
  $n$ ``inverse temperatures" $\beta_1, ..., \beta_n$, which we
  combine into one vector quantity $\beta$ lying in the dual space $\RR^{n*}$.
  
  \vskip .2cm
  
  For simplicity we assume that $E$ defines an embedding of $A$ into $\RR^n$.
  We can then think of $A$ as being a subset of $\RR^n$ to begin with,
  and sometimes drop $E$ from the notation,
  thinking of it as just the inclusion map.
  With these conventions, we write the partition function and the Gibbs distribution
  \be\label{eq:gibbs-dist-vect}
  Z(\beta)\,\,=\,\,\sum_{\omega\in A} e^{-(\beta, \omega)}, \quad 
  p_\omega(\beta) = {e^{-(\beta, \omega)}\over Z(\beta)}, 
  \quad \beta\in\RR^{n*}.
  \ee
  As in \eqref{eq:mean-value-O}, the Gibbs distribution can be used to define
  the mean value of any observable $O$ on $A$. In particular, taking for $O$
  the vector valued function (embedding) $E: A\to\RR^n$, we have the 
  {\em mean energy map}
  \be
  \la E\ra : \beta\,\,\longmapsto \la E\ra(\beta) \,\,=\,\,{\sum_{\omega\in A} \omega\cdot e^{-(\beta, \omega)}\over
  \sum_{\omega'\in A} e^{-(\beta,\omega)}}\,\,=\,\,-\nabla_\beta \log Z(\beta).
  \ee
  Here $\nabla_\beta$ means the vector of gradient with respect to $\beta$, i.e., the
  differential of a function considered as a vector in the dual space. Thus $\la E\ra: \RR^{n*}\to\RR^n$. 
  Let $Q\subset \RR^n$ be the convex hull of $A$, and $Q^\circ$ be the interior of $Q$.
  Since $(p_\omega(\beta))_{\omega\in A}$ is a probability distribution on $A$ with all components
  nonzero, we see that $\la E\ra$ maps $\RR^{n*}$ into $Q^\circ$. We can now generalize the
  thermodynamic formalism of the previous section as follows.
  
  \begin{prop}\label{prop:diffeo-2}
  (a) The map $\la E\ra: \RR^{n*}\to Q^\circ$ is a diffeomorphism.
  
  (b) For any $\overline E\in Q^\circ$ let  $P_{\overline E}$ be the set of probability distributions
  $p\in\Delta^A$ satisfying the constraints
  $$\la E\ra_p \,\,:=\,\,\sum_{\omega\in A} p_\omega\cdot\omega \,\,=\,\,\overline E.
  $$
  Let $\beta(\overline E) \in\RR^{n*}$ be unique such that 
  $\la E\ra(\beta)=\overline E$. Then the Gibbs distribution $p(\beta) = (p_\omega(\beta))$
  defined above, has maximal entropy among all the distributions from $P_{\overline E}$.
  
  (c) Let $S(\overline E)$ be the entropy of the distribution $p(\beta(\overline E))$.
  Then $\nabla_{\overline E} S(\overline E) = \beta(\overline E)$. 
  
  (d) The functions $-\log Z(\beta)$ on $\RR^{n*}$ and $S(\overline E)$ on $Q^\circ\subset \RR^n$
  are concave and are the Legendre transforms of each other. 
   \end{prop} 
   
   Note that  part (b) shows that the ``vector" Gibbs distribution \eqref{eq:gibbs-dist-vect} has the same
   thermodynamic significance as the more standard one \eqref{eq:gibbs-one}.

   \vskip .2cm
   
  The proof of the proposition will be given later in this section.
  
  \paragraph{B. Example: toric varieties and the moment map.} Assume that $A$ lies in $\ZZ^n\subset\RR^n$. 
  The exponential $e^{-(\beta, \omega)},\omega\in A$, then becomes a Laurent monomial
  $z^\omega = \prod z_i^{\omega_i}$ in the variables $z_i = e^{-\beta_i}$. Real values of $\beta$
  correspond to $z\in \RR^n_+$, where $\RR_+$ is the set of positive real numbers.
  
  \vskip .2cm

The monomial $z^\omega$ makes sense for any $z\in(\CC^*)^n$. Consider the complex vector
space $\CC^A$ with basis $e_\omega, \omega\in A$ and let $\PP^A$ be its projectivization.
A vector of $\CC^A$ is thus a tuple $(a_\omega)_{\omega\in A}$.
The torus $(\CC^*)^n$ acts on $\CC^A$ and $\PP^A$ by
$$z\cdot  e_\omega\,\,=\,\, z^\omega e_\omega.$$
In particular, we consider the orbit $X_A^\circ\subset \PP^A$ of the point 
represented by $\1=(1)_{\omega\in A}\in\CC^A$ and let $X_A\subset \PP^A$ be
the projective toric variety defined as the closure of $X_A^\circ$. 

\vskip .2cm

Assume for simplicity that $A$ generates $\ZZ^n$ as an affine lattice, i.e., there is no
smaller integer affine sublattice in $\ZZ^n$ containing $A$. Then the action of $\CC^n$ on 
$\PP^A$ is faithful, in particular, the action map
$$z\,\,\longmapsto \,\, z\cdot \1 \,\,=\,\,(z^\omega)_{\omega\in A}
$$
identifies $(\CC^*)^n$ with $X_A^\circ$. Let $X_A^+\subset X_A^\circ$ be the image of
$\RR_+^n\subset (\CC^*)^n$. This image is known as the {\em positive part}
of the toric variety $X_A$. Clearly, $X_A^+$ consists of the points of the form
$x(\beta) = (e^{-(\beta,\omega)})_{\omega\in A}$ for all $\beta\in\RR^{n*}$. 

\vskip .2cm

The action of the compact part $(S^1)^n$ of the torus $(\CC^*)^n$ on the projective space
$\PP^A$ preserves the standard Fubini-Study K\"ahler metric and gives rise to the 
{\em moment map}
\be
\mu_\PP: \PP^A \to \RR^n, \quad (a_\omega)_{\omega\in A} \,\,\longmapsto \,\, {\sum_{\omega\in A} 
\omega\cdot \|a_\omega\|^2\over \sum_{\omega\in A} \|a_\omega\|^2},
\ee
see, e.g., \cite{GKZ}.  The image of this map is the polytope $Q=\operatorname{Conv}(A)$. 
Let  $ \mu_X^+$ be the restriction of $\mu_P$ to  $X_A^+$. Using the above
parametrization of $X_A^+$ by the $x(\beta)$, we write $\mu_X^+$ as a map
from $\RR^{n*}$ to $Q$, and find that 
\be
\mu_A^+(\beta) \,\,=\,\, {\sum_{\omega\in A} 
\omega\cdot e^{-(2\beta, \omega)} \over \sum_{\omega\in A} e^{-(2\beta,\omega)}}\,\,=\,\,
\la E\ra (2\beta)
\ee
is nothing but the mean energy  of the twice scaled $\beta$, with respect to the Gibbs distribution
 \eqref{eq:gibbs-dist-vect}.
Proposition \ref{prop:diffeo-2}(a) reduces then to the well known fact about
toric varieties: that the moment map defines a diffeomorphism from the positive part
to the interior of the defining polytope, see \cite{atiyah}, \cite{GKZ}.

\paragraph{C. Direct and inverse images of concave functions.}
To give a natural proof of Proposition  \ref{prop:diffeo-2}, we start with some general remarks.
By a {\em convex body} we will mean a convex subset $P$ of  
 some finite-dimensional affine space $V$ over $\RR$. For such $P$
   we denote by $\Conc(P)$ the space (semigroup)
 of concave functions $f: P\to\RR$ which are proper, i.e., such that each
 level set $f^{-1}(c)$ is compact.
 Any such function achieves a maximum on $P$.

 By an {\em admissible  embedding}
 of convex bodies $i: P'\to P$ we mean an injective map induced by an
 affine embedding of ambient affine spaces $V'\hookrightarrow V$, 
 so that $P'=V'\cap P$. In this case for any $f\in\Conc(P)$ we have the
 {\em inverse image} (restriction) $i^*f=f|_{P'}$ which again lies in $\Conc(P')$.

 Similarly, by an {\em admissible surjection}
 of convex bodies $j: P\to P''$ we mean a surjective map induced by
 an affine surjection $J: A'\to A''$  of ambient affine spaces. In this case for
 any $f\in \Conc(P)$ we have the {\em direct image} which is the function
 $j_*f$ on $P''$ defined by 
 \be
(j_* f)(p'') \,\,=\,\,\max_{j(p)=p''} f(p).
\ee

\begin{ex}
 One can take $P=V$ to be  a finite-dimensional vector space over $\RR$ and
 $f$ to be a negative definite quadratic form on $V$. Then, for any linear
 surjection $j: V\to V''$,  the direct image $j_*f$ is a negative definite quadratic
 form on $V''$.  The integration, along the fibers of $j$,  of the Gaussian
 function $e^{f(p)}$ on $V$ gives, up to a constant, the Gaussian function
 $e^{(j_*f)(p'')}$ on $V''$. 
 
 For a general $f\in \Conc(P)$ and an admissible surjection $j: P\to P''$
 the function $e^{j_*f}$ is the leading term, as $h\to 0$, of the  function on $P''$ obtained by
 integrating $e^{f(p)/h}$ along the fibers of $j$.
  \end{ex}
 The following is then elementary.
 
 \begin{prop} (a) The function $j_*f$   belongs to $\Conc(P'')$.
 
 (b) (Base change)  Let
$$
\xymatrix{P_2 \ar[r]^{i} \ar[d]_{j_2} &P_1\ar[d]^{j_1}\cr
P'_2 \ar[r]^{i'}&P'_1
}
$$
be a Cartesian square of convex bodies, such that $i, i'$
are admissible embeddings and $j_1, j_2$ are 
admissible surjections. Then for any  $f\in\Conc(P_1)$ we have the equality
 $(j_1)_* i^*f = (i')^* (j_2)_* f$  of concave functions on $P'_2$. 
 \qed 
 \end{prop}
 
 \paragraph{D. Proof of Proposition \ref{prop:diffeo-2}.} Consider the admissible
 surjection of convex bodies
 \be
 \pi: \Delta^A\lra Q, \quad p\,\,\longmapsto \,\, \la E\ra_p\,\, =  \,\, \sum_{\omega\in A}\,  p_\omega\cdot \omega.
 \ee
 Fix $\overline E\in Q^\circ$. The fiber $\pi^{-1}(\overline E)$ is the set $P_{\overline E}$ of part (b) of the proposition. 
  Consider the entropy function $S\in\Conc(\Delta^A)$ defined by 
 \eqref{eq:entropy-def}. It is strictly concave, so the restriction of $S$ to 
 $\pi^{-1}(\overline E)$ achieves maximum at a unique interior point; denote this
 point  $\widetilde p(\overline E)$. This point is, furthermore,  the unique critical point  of $S$
 on  $\pi^{-1}(\overline E)$. Consider also the direct image function $\pi_*S\in\Conc(Q)$.

  The location of the critical point can be found by the Lagrange
 multiplier method for $n$ constraints: the differential of $S$ at  $\widetilde p(\overline E)$
 should be a linear combination of the differentials of the individual scalar constraints,
  i.e., to have the form
 $(\lambda, d\pi)$ for some $\lambda\in \RR^{n*}$.     Further,  we have the $n$-constraint
 interpretation of the  Lagrange multipliers as  minus the partial derivatives of the maximal value
 with respect to the constraints, see again \cite{econ-book}. 
This means that $\lambda = \lambda(\overline E)$ is equal to the gradient
 (differential) of the strictly convex function
 of $\pi_*S$ at the point  $\overline E$.

 Next, look at the Gibbs distribution $p(\lambda(\overline E))$. We see that 
 $p(\lambda(\overline E))$ is  a critical point of $S$ on $\pi^{-1}(\overline E)$,
 so it is equal to $\widetilde p(\overline E)$. This implies that $\lambda(\overline E)=\beta(\overline E)$
 is the inverse to the map $\beta\mapsto \la E\ra(\beta)$ which is therefore a diffeomorphism,
 thus proving part (a) of the proposition.  In particular, the function $S(\overline E)$ of part (c)
 is the same as $\pi_*S$. 
 Since $p(\lambda(\overline E))=\widetilde p(\overline E)$,
 this implies part (b). Part (c) follows since $\lambda(\overline E) = \nabla_{\overline E} (\pi_*S)$
 by definition. Finally, the Legendre transform relation between the functions $-\log Z(\beta)$ 
 and $S(\overline E)=\pi_*S$ in part (d) is equivalent to the fact
 that their gradients define mutually inverse diffeomorphisms, as we have shown that
 $\la E\ra = \nabla_\beta (-\log Z)$ is inverse to $\nabla_{\overline E}(\pi_*S)$. \qed

\vskip 2cm

 Department of Mathematics, Yale University, 10 Hillhouse Avenue, New Haven CT 06520 USA, 
 \texttt{ mikhail.kapranov@yale.edu}

\ed

\end{document}